\documentclass[12pt]{article}
\usepackage{amssymb,amsmath}
\usepackage{cases}
\usepackage{amsfonts}
\usepackage{color,xcolor}
\usepackage[left=2.0cm,right=2.0cm,top=2.0cm,bottom=2.0cm]{geometry}
\usepackage[colorlinks,citecolor=blue,urlcolor=blue]{hyperref}

\newtheorem{theorem}{Theorem}[section]

\newtheorem{lemma}{Lemma}[section]
\newtheorem{proposition}{Proposition}[section]

\newtheorem{definition}{Definition}[section]
\newtheorem{remark}{Remark}[section]

\newcommand{\bal}{\begin{align}}
\newcommand{\bbal}{\begin{align*}}
\newcommand{\beq}{\begin{equation}}
\newcommand{\eeq}{\end{equation}}
\newcommand{\bca}{\begin{cases}}
\newcommand{\eca}{\end{cases}}

\newcommand{\dd}{\mathrm{d}}

\newcommand{\R}{\mathbb{R}}

\newcommand{\bi}{\Big}

\begin{document}

\title{Non-uniform continuous dependence on initial data  for a two-component Novikov system  in Besov space}

\author{Xing Wu\thanks{Corresponding author. ny2008wx@163.com (Xing Wu)},  Jie Cao\\
\small \it College of Information and Management Science,
Henan Agricultural University,\\
\small Zhengzhou, Henan, 450002, China}

\date{}

\maketitle\noindent{\hrulefill}

{\bf Abstract:} In this paper, we show that the solution map of the two-component Novikov system  is not uniformly continuous on the initial data in Besov spaces $B_{p, r}^{s-1}(\mathbb{R})\times B_{p, r}^s(\mathbb{R})$ with $s>\max\{1+\frac{1}{p}, \frac{3}{2}\}$, $1\leq p< \infty$,  $1\leq r<\infty$. Our result covers and extends the previous non-uniform continuity in Sobolev spaces $H^{s-1}(\mathbb{R})\times H^s(\mathbb{R})$ for $s>\frac{5}{2}$ (J. Math. Phys., 2017) to Besov spaces.

{\bf Keywords:} Non-uniform dependence; two component Novikov system; Besov spaces

{\bf MSC (2010):} 35B30; 35G25; 35Q53
\vskip0mm\noindent{\hrulefill}

\section{Introduction}\label{sec1}
In this paper, we are concerned with the Cauchy problem for the following  two-component Novikov system on $\mathbb{R}$:
\begin{eqnarray}\label{eq1}
        \left\{\begin{array}{ll}
         \rho_t=\rho_xu^2+\rho uu_x,\\
          m_t=3u_xum+u^2m_x-\rho(u\rho)_x,\\
          m=u-u_{xx},\\
         \rho(0, x)=\rho_0, u(0, x)=u_0. \end{array}\right.
        \end{eqnarray}
The system in (\ref{eq1}) was proposed by Popowicz in \cite{Popowicz 2015} as  the two-component generalization of the Novikov equation and can be rewritten in the Hamiltonian form (see \cite{Popowicz 2015} for details).

If we take $\rho=0$, the system in (\ref{eq1}) reduces to
\begin{equation}\label{eq2}
 m_t=3u_xum+u^2m_x, \qquad m=u-u_{xx},
\end{equation}
which is nothing but the famous Novikov equation derived in \cite{Novikov 2009} as a new integrable equation with cubic nonlinearities. It is shown in \cite{Home 2008} that the Novikov equation admits peakon solutions, has a Lax pair in matrix form and a bi-Hamiltonian structure. Furthermore, it has infinitely many conserved quantities.

The local well-posedness and dependence on initial data for the Novikov equation (\ref{eq2}) in Sobolev spaces and Besov spaces were established in \cite{Ni 2011, Himonas 2012, Wu 2012, Yan 2012, Wu 2013}. After the non-uniform dependence  for some dispersive equations was studied by Kenig et al. \cite{Kenig 2001}, the issue of non-uniform continuity of solutions on initial data has attracted much more attention, such as on classical Camassa-Holm equation \cite{Himonas 2005, Himonas 2009, Himonas 2010, Li 2020} and on famous Degasperis-Procesi equation \cite{Christov 2009, Himonas 2011, Wu 2020}.  The first result of nonuniform dependence for the  Novikov equation (\ref{eq2}) was proved  by Himonas and Holliman  \cite{ Himonas 2012} in Sobolev space $H^s(\mathbb{R})$, $s>\frac{3}{2}$. Recently, Li, Li and Zhu \cite {1Li 2020} have improved the result in \cite{ Himonas 2012} to Besov spaces $B_{p, r}^s(\mathbb{R})$ with  $s>\max\{1+\frac{1}{p}, \frac{3}{2}\}$, $1\leq p\leq \infty$, $1\leq r<\infty$.

Compared with the rich research on the Novikov equation, there are few  mathematical studies on the two-component Novikov system. The local well-posedness for system (\ref{eq1}) in Besov space was first investigated by Luo and Yin \cite{Luo 2015} with initial data $(\rho_0, u_0)\in B_{p, r}^{s-1}(\mathbb{R})\times B_{p, r}^s(\mathbb{R})$ for $s>\max\{1+\frac{1}{p}, \frac{3}{2}\}$, $1\leq p, r\leq \infty$ or for the critical index $(s, p, r)=(\frac{3}{2}, 2, 1)$. They also showed that  the solution map is H\"{o}lder continuous from $B_{p, r}^{s-1}\times B_{p, r}^s$ to $B_{p, r}^{s'-1}\times B_{p, r}^{s'}$
with $s-1\leq s' < s$. Actually, the continuity of the solution map on the initial data (i.e. $s'=s$) can be achieved by using the theory established in \cite{Li 2016}. Moreover, two blow-up criteria for the system were presented by making use of the conservation laws. Wang and Fu \cite{Wang 2017} further proved that the solution map is not uniformly continuous in Sobolev spaces $H^{s-1}(\mathbb{R})\times H^s(\mathbb{R})$ with $s>\frac{5}{2}$.

In the present paper, motivated by \cite{Li 2020},  we aim at showing that  the solution map of (\ref{eq1}) is not uniformly continuous depending on the initial data in Besov spaces $B_{p, r}^{s-1}\times B_{p, r}^s$ with $s>\max\{1+\frac{1}{p}, \frac{3}{2}\}$, $1\leq p < \infty$, $1\leq r< \infty$.

For studying the non-uniform continuity of the two-component Novikov system, it is more convenient to express (\ref{eq1}) in the following equivalent nonlocal form
\begin{eqnarray}\label{eq3}
        \left\{\begin{array}{ll}
        \rho_t=u^2\rho_x+\rho uu_x,\\
         u_t=u^2u_x+\mathcal{P}(u)+\mathcal{R}(u, \rho),\\
          \rho(0, x)=\rho_0, u(0, x)=u_0,\end{array}\right.
        \end{eqnarray}
where $\mathcal{P}(u)=\mathcal{P}_1(u)+\mathcal{P}_2(u)+\mathcal{P}_3(u),$ $\mathcal{R}(u, \rho)=\mathcal{R}_1(u, \rho)+\mathcal{R}_2(u, \rho) $ and
\begin{eqnarray*}
&\;&\mathcal{P}_1(u)=\partial_x(1-\partial_x^2)^{-1}(u^3), \quad \mathcal{P}_2(u)=\frac{3}{2}\partial_x(1-\partial_x^2)^{-1}(uu_x^2), \quad \mathcal{P}_3(u)=\frac{1}{2}(1-\partial_x^2)^{-1}(u_x^3),\\
&\;&\mathcal{R}_1(u, \rho)=-\frac{1}{2}\partial_x(1-\partial_x^2)^{-1}(u\rho^2),\quad
\mathcal{R}_2(u, \rho)=-\frac{1}{2}(1-\partial_x^2)^{-1}(u_x\rho^2).
\end{eqnarray*}

 Our main result is stated as follows.
\begin{theorem}\label{the1.1} Let $s>\max\{1+\frac{1}{p}, \frac{3}{2}\}$, $1\leq p < \infty$,  $1\leq r<\infty$. The solution map $(\rho_0, u_0)\rightarrow (\rho(t), u(t))$ of the initial value problem (\ref{eq3}) is not
uniformly continuous from any bounded subset of  $B_{p, r}^{s-1}(\mathbb{R})\times B_{p, r}^s(\mathbb{R})$ into $\mathcal{C}([0, T];  B_{p, r}^{s-1}(\mathbb{R})\times B_{p, r}^s(\mathbb{R}))$. More precisely,
there exist two sequences of solutions $(\rho_n, u_n)$ and $(\tilde{\rho}_n, v_n)$ such that the corresponding initial data satisfy
\begin{eqnarray*}
        \|\rho_{0, n}, \tilde{\rho}_{0, n}\|_{B_{p, r}^{s-1}}+\|u_{0, n}, v_{0, n}\|_{B_{p, r}^s}\lesssim 1,\\
\end{eqnarray*}
       \mbox{and}
   \begin{eqnarray*}
        \lim_{n\rightarrow \infty}( \|\tilde{\rho}_{0, n}-\rho_{0, n}\|_{B_{p, r}^{s-1}}+\|v_{0, n}-u_{0, n}\|_{B_{p, r}^s})=0,
        \end{eqnarray*}
but
\begin{eqnarray*}
  \liminf_{n\rightarrow \infty} \|\tilde{\rho}_n-\rho_n\|_{B_{p, r}^{s-1}}\gtrsim t, \quad    \liminf_{n\rightarrow \infty} \|v_n-u_n\|_{B_{p, r}^s}\gtrsim t, \quad t\in [0, T_0],
        \end{eqnarray*}
with small positive time $T_0$ for $T_0\leq T$.
\end{theorem}
\begin{remark}\label{rem1}
Since $B_{2, 2}^s=H^s$, our result covers and extends the previous non-uniform continuity of solutions on initial data in Sobolev spaces $H^{s-1}(\mathbb{R})\times H^s(\mathbb{R})$ for $s>\frac{5}{2}$ \cite{Wang 2017} to Besov spaces.
\end{remark}

{\bf Notations}:  Given a Banach space $X$, we denote the norm of a function on $X$ by $\|\|_{X}$, and \begin{eqnarray*}
\|\cdot\|_{L_T^\infty(X)}=\sup_{0\leq t\leq T}\|\cdot\|_{X}.
\end{eqnarray*}
For $\mathbf{f}=(f_1, f_2,...,f_n)\in X$,
\begin{eqnarray*}
\|\mathbf{f}\|_{X}^2=\|f_1\|_{X}^2+\|f_2\|_{X}^2+...+\|f_n\|_{X}^2.
\end{eqnarray*}
The symbol
$A\lesssim B$ means that there is a uniform positive constant $C$ independent of $A$ and $B$ such that $A\leq CB$.

\section{Littlewood-Paley analysis}\label{sec2}
\setcounter{equation}{0}
In this section, we will review the definition of Littlewood-Paley decomposition and nonhomogeneous Besov space, and then list some useful properties. For more details, the readers can refer to \cite{Bahouri 2011}.

There exists a couple of smooth functions $(\chi,\varphi)$ valued in $[0,1]$, such that $\chi$ is supported in the ball $\mathcal{B}\triangleq \{\xi\in\mathbb{R}^d:|\xi|\leq \frac 4 3\}$, $\varphi$ is supported in the ring $\mathcal{C}\triangleq \{\xi\in\mathbb{R}^d:\frac 3 4\leq|\xi|\leq \frac 8 3\}$. Moreover,
$$\forall\,\, \xi\in\mathbb{R}^d,\,\, \chi(\xi)+{\sum\limits_{j\geq0}\varphi(2^{-j}\xi)}=1,$$
$$\forall\,\, \xi\in\mathbb{R}^d\setminus\{0\},\,\, {\sum\limits_{j\in \mathbb{Z}}\varphi(2^{-j}\xi)}=1,$$
$$|j-j'|\geq 2\Rightarrow\textrm{Supp}\,\ \varphi(2^{-j}\cdot)\cap \textrm{Supp}\,\, \varphi(2^{-j'}\cdot)=\emptyset,$$
$$j\geq 1\Rightarrow\textrm{Supp}\,\, \chi(\cdot)\cap \textrm{Supp}\,\, \varphi(2^{-j}\cdot)=\emptyset.$$
Then, we can define the nonhomogeneous dyadic blocks $\Delta_j$ as follows:
$$\Delta_j{u}= 0,\,\, \text{if}\,\, j\leq -2,\quad
\Delta_{-1}{u}= \chi(D)u=\mathcal{F}^{-1}(\chi \mathcal{F}u),$$
$$\Delta_j{u}= \varphi(2^{-j}D)u=\mathcal{F}^{-1}(\varphi(2^{-j}\cdot)\mathcal{F}u),\,\, \text{if} \,\, j\geq 0.$$

\begin{definition}[\cite{Bahouri 2011}]\label{de2.1}
Let $s\in\mathbb{R}$ and $1\leq p,r\leq\infty$. The nonhomogeneous Besov space $B^s_{p,r}(\mathbb{R}^d)$ consists of all tempered distribution $u$ such that
\begin{align*}
||u||_{B^s_{p,r}(\mathbb{R}^d)}\triangleq \Big|\Big|(2^{js}||\Delta_j{u}||_{L^p(\mathbb{R}^d)})_{j\in \mathbb{Z}}\Big|\Big|_{\ell^r(\mathbb{Z})}<\infty.
\end{align*}
\end{definition}

In the following, we list some basic lemmas and properties about Besov space which will be frequently used in proving our main result.

\begin{lemma}(\cite{Bahouri 2011})\label{lem2.1}
 (1) Algebraic properties: $\forall s>0,$ $B_{p, r}^s(\mathbb{R}^d)$ $\cap$ $L^\infty(\mathbb{R}^d)$ is a Banach algebra. $B_{p, r}^s(\mathbb{R}^d)$ is a Banach algebra $\Leftrightarrow B_{p, r}^s(\mathbb{R}^d)\hookrightarrow L^\infty(\mathbb{R}^d)\Leftrightarrow s>\frac{d}{p}$ or $s=\frac{d}{p},$ $r=1$.\\
 (2) For any $s>0$ and $1\leq p,r\leq\infty$, there exists a positive constant $C=C(d,s,p,r)$ such that
$$\|uv\|_{B^s_{p,r}(\mathbb{R}^d)}\leq C\Big(\|u\|_{L^{\infty}(\mathbb{R}^d)}\|v\|_{B^s_{p,r}(\mathbb{R}^d)}+\|v\|_{L^{\infty}(\mathbb{R}^d)}\|u\|_{B^s_{p,r}(\mathbb{R}^d)}\Big).$$
(3) Let $m\in \mathbb{R}$ and $f$ be an $S^m-$ multiplier (i.e., $f: \mathbb{R}^d\rightarrow \mathbb{R}$ is smooth and satisfies that $\forall \alpha\in \mathbb{N}^d$, there exists a constant $\mathcal{C}_\alpha$ such that $|\partial^\alpha f(\xi)|\leq \mathcal{C}_\alpha(1+|\xi|)^{m-|\alpha|}$ for all $\xi \in \mathbb{R}^d$). Then the operator $f(D)$ is continuous from $B_{p, r}^s(\mathbb{R}^d)$ to $B_{p, r}^{s-m}(\mathbb{R}^d)$.\\
(4) Let  $1\leq p, r\leq \infty$ and $s>\max\{1+\frac{1}{p}, \frac{3}{2}\}$. Then  we have
$$\|uv\|_{B_{p, r}^{s-2}(\mathbb{R})}\leq C\|u\|_{B_{p, r}^{s-2}(\mathbb{R})}\|v\|_{B_{p, r}^{s-1}(\mathbb{R})}.$$
Hence, for the operators $\mathcal{P}$ and $\mathcal{R}$ in (\ref{eq3}), we have
\begin{align*}
  \|\mathcal{P}(u)-\mathcal{P}(v)\|_{B_{p, r}^{s-1}(\mathbb{R})}&\lesssim \|u-v\|_{B_{p, r}^{s-1}(\mathbb{R})}\|u, v\|^2_{B_{p, r}^s(\mathbb{R})},\\
 \|\mathcal{R}(\tilde{u}, \tilde{v})-\mathcal{R}(u, v)\|_{B_{p, r}^{s-1}(\mathbb{R})}&\lesssim \|\tilde{u}-u\|_{B_{p, r}^{s-1}(\mathbb{R})}\|\tilde{v}\|^2_{B_{p, r}^{s-1}(\mathbb{R})}+\|\tilde{v}-v\|_{B_{p, r}^{s-2}(\mathbb{R})}(\|u\|^2_{B_{p, r}^s(\mathbb{R})}+\|\tilde{v}, v\|^2_{B_{p, r}^{s-1}(\mathbb{R})}).
\end{align*}
\end{lemma}

\begin{lemma}\label{lem2.2}(\cite{Bahouri 2011, Li 2017})
Let $1\leq p,r\leq \infty$. Assume that
\begin{eqnarray*}
\sigma> -d \min\{\frac{1}{p}, 1-\frac{1}{p}\} \quad \mathrm{or}\quad \sigma> -1-d \min\{\frac{1}{p}, 1-\frac{1}{p}\}\quad \mathrm{if} \quad \mathrm{div\,} v=0.
\end{eqnarray*}
There exists a constant $C=C(d,p,r,\sigma)$ such that for any solution to the
following linear transport equation:
\begin{equation*}
\partial_t f+v\cdot\nabla f=g,\qquad
f|_{t=0} =f_0,
\end{equation*}
the following statements hold:
\begin{align*}
\sup_{s\in [0,t]}\|f(s)\|_{B^{\sigma}_{p,r}}\leq Ce^{CV_{p}(v,t)}\Big(\|f_0\|_{B^\sigma_{p,r}}
+\int^t_0\|g(\tau)\|_{B^{\sigma}_{p,r}}\dd \tau\Big),
\end{align*}
with
\begin{align*}
V_{p}(v,t)=
\begin{cases}
\int_0^t \|\nabla v(s)\|_{B^{\sigma-1}_{p,r}}\dd s, &\quad \mathrm{if} \;\sigma>1+\frac{d}{p}\ \mathrm{or}\ \{\sigma=1+\frac{d}{p} \mbox{ and } r=1\},\\
\int_0^t \|\nabla v(s)\|_{B^{\sigma}_{p,r}}\dd s,&\quad\mathrm{if} \; \sigma=1+\frac{d}{p} \quad \mathrm{and} \quad r>1,\\
\int_0^t \|\nabla v(s)\|_{B^{\frac{d}{p}}_{p,\infty}\cap L^\infty}\dd s,&\quad\mathrm{if} \; \sigma<1+\frac{d}{p}.
\end{cases}
\end{align*}
\end{lemma}

\section{Non-uniform continuous dependence}\label{sec3}
\setcounter{equation}{0}
In this section, we will give the proof of our main theorem.

Let $\hat{\phi}\in \mathcal{C}^\infty_0(\mathbb{R})$ be an even, real-valued and non-negative funtion on $\R$ and satisfy
\begin{numcases}{\hat{\phi}(x)=}
1, &if $|x|\leq \frac{1}{4}$,\nonumber\\
0, &if $|x|\geq \frac{1}{2}$.\nonumber
\end{numcases}
Define the high frequency function $f_n$ and the low frequency function $g_n$ by
$$f_n=2^{-ns}\phi(x)\sin \bi(\frac{17}{12}2^nx\bi), \qquad g_n=2^{-\frac n2}\phi(x), \quad n\gg1.$$
It has been showed in \cite{Li 2020} that $\|f_n\|_{B_{p, r}^\sigma}\lesssim 2^{n(\sigma-s)}$.

It is easy to verify that $(2^nf_n, f_n)$ lies in $B_{p, r}^{s-1}(\mathbb{R})\times B_{p, r}^{s}(\mathbb{R})$, so does $(2^nf_n+g_n, f_n+g_n)$.\\
 Let
$$(\rho_{0, n}, u_{0, n})=(2^nf_n, f_n), \quad (\tilde{\rho}_{0, n}, v_{0, n})=(2^nf_n+g_n, f_n+g_n),$$
we have
\begin{align}
\|\tilde{\rho}_{0, n}\|_{B_{p, r}^{s-2}}&\lesssim 2^{-\frac n2},\quad  \|\tilde{\rho}_{0, n}\|_{B_{p, r}^{s+\sigma}}\lesssim 2^{n(\sigma+1)}, \; \sigma\geq -\frac 32, \label{eq3.1}\\
 \|v_{0, n}\|_{B_{p, r}^{s-1}}&\lesssim 2^{-\frac n2},\quad  \|v_{0, n}\|_{B_{p, r}^{s+\sigma}}\lesssim 2^{n\sigma}, \; \sigma\geq -\frac 12. \label{eq3.2}
\end{align}

Consider the system (\ref{eq3}) with initial data $(\rho_{0, n}, u_{0, n})$ and $(\tilde{\rho}_{0, n}, v_{0, n})$, respectively.  According to the local well-posedness result in \cite{Luo 2015}, there exists corresponding solution $(\rho_n, u_n)$, $(\tilde{\rho}_n, v_n)$ belonging to
$\mathcal{C}([0, T];  B_{p, r}^s)$ and has common lifespan $T\thickapprox 1$. Moreover, by Lemma \ref{lem2.1}-\ref{lem2.2}, there holds
 \begin{eqnarray}
     \|\rho_n\|_{L_T^\infty(B_{p, r}^{s+k-1})}+\|u_n\|_{L_T^\infty(B_{p, r}^{s+k})}\lesssim  \|\rho_{0, n}\|_{B_{p, r}^{s+k-1}}+\|u_{0, n}\|_{B_{p, r}^{s+k}}\lesssim2^{nk},\quad k\geq -1,\label{eq3.3}\\
      \|\tilde{\rho}_n\|_{L_T^\infty(B_{p, r}^{s+l-1})}+\|v_n\|_{L_T^\infty(B_{p, r}^{s+l})}\lesssim  \|\tilde{\rho}_{0, n}\|_{B_{p, r}^{s+l-1}}+\|v_{0, n}\|_{B_{p, r}^{s+l}}\lesssim\left\{\begin{array}{ll}
       2^{nl},\quad\;\; l\geq -\frac{1}{2}, \\
        2^{-\frac n2}, \quad l=-1.
        \end{array}\right.\label{eq3.4}
        \end{eqnarray}

In the following, we shall firstly show that for the selected high  frequency initial data $(\rho_{0, n}, u_{0, n})$, the corresponding solution $(\rho_n, u_n)$ can be approximated by the initial data. More precisely, that is
\begin{proposition}\label{pro1}
 Under the assumptions of Theorem \ref{the1.1}, we have
      \begin{eqnarray}
      \|\rho_n-\rho_{0, n}\|_{L_T^\infty(B_{p, r}^{s-1})}+\|u_n-u_{0, n}\|_{L_T^\infty(B_{p, r}^s)}\lesssim 2^{-\frac{n}{2}(s-1)}.\label{eq3.5}
        \end{eqnarray}
\end{proposition}
\noindent{\bf Proof} \;
Denote $\varrho=\rho_n-\rho_{0, n}$, $\epsilon=u_n-u_{0, n}$, then we can derive from (\ref{eq3}) that $(\varrho, \epsilon)$ satisfies
\begin{eqnarray}\label{eq3.6}
        \left\{\begin{array}{ll}
        \varrho_t-(u_n)^2\partial_x\varrho=[(u_n)^2-(u_{0, n})^2]\partial_x\rho_{0, n}+(u_{0, n})^2\partial_x\rho_{0, n}+\rho_nu_n\partial_xu_n,\\
         \epsilon_t-(u_n)^2\partial_x\epsilon=[(u_n)^2-(u_{0, n})^2]\partial_xu_{0, n}+[\mathcal{P}(u_n)-\mathcal{P}(u_{0, n})]+[\mathcal{R}(u_n, \rho_n)-\mathcal{R}(u_{0, n}, \rho_{0, n})]\\
         \qquad\qquad\qquad\;\;\;+\mathcal{P}(u_{0, n})+\mathcal{R}(u_{0, n}, \rho_{0, n})+(u_{0, n})^2\partial_xu_{0, n},\\
          \varrho(0, x)=0, \epsilon(0, x)=0,\end{array}\right.
        \end{eqnarray}

Applying Lemma \ref{lem2.2} yields
\begin{align}
     \|\varrho\|_{B_{p, r}^{s-2}}&\lesssim \int_0^t\|\partial_x(u_n)^2\|_{B_{p, r}^{s-1}}\|\varrho\|_{B_{p, r}^{s-2}}d\tau+\int_0^t \|[(u_n)^2-(u_{0, n})^2]\partial_x\rho_{0, n}\|_{B_{p, r}^{s-2}}d\tau \nonumber\\
    &\quad +\int_0^t \|\rho_nu_n\partial_xu_n-\rho_{0, n}u_{0, n}\partial_xu_{0, n}\|_{B_{p, r}^{s-2}}d\tau\nonumber\\
    &\quad  +t\|(u_{0, n})^2\partial_x\rho_{0, n}, \rho_{0, n}u_{0, n}\partial_xu_{0, n}\|_{B_{p, r}^{s-2}},\label{eq3.7}\\
\|\epsilon\|_{B_{p, r}^{s-1}}&\lesssim \int_0^t\|\partial_x(u_n)^2\|_{B_{p, r}^{s-1}}\|\epsilon\|_{B_{p, r}^{s-1}}d\tau+\int_0^t \|[(u_n)^2-(u_{0, n})^2]\partial_xu_{0, n}\|_{B_{p, r}^{s-1}}d\tau \nonumber\\
    &\quad +\int_0^t \|[\mathcal{P}(u_n)-\mathcal{P}(u_{0, n})], [\mathcal{R}(u_n, \rho_n)-\mathcal{R}(u_{0, n}, \rho_{0, n})]\|_{B_{p, r}^{s-1}}d\tau \nonumber\\
     &\quad+t\|\mathcal{P}(u_{0, n}), \mathcal{R}(u_{0, n}, \rho_{0, n}), (u_{0, n})^2\partial_xu_{0, n}\|_{B_{p, r}^{s-1}},\label{eq3.8}
\end{align}

Using Lemma \ref{lem2.1} and  the fact that $B_{p, r}^{s-1}(\mathbb{R})$ is a  Banach algebra when $s>\max\{1+\frac{1}{p}, \frac{3}{2}\}$, we have
\begin{align*}
      \|[(u_n)^2-(u_{0, n})^2]\partial_x\rho_{0, n}\|_{B_{p, r}^{s-2}}&\lesssim  \|\partial_x\rho_{0, n}\|_{B_{p, r}^{s-2}}\|u_n-u_{0, n}\|_{B_{p, r}^{s-1}}\|u_n-u_{0, n}\|_{B_{p, r}^{s-1}}\\
      &\lesssim \|u_n-u_{0, n}\|_{B_{p, r}^{s-1}}(\|\rho_{0, n}, u_{0, n}, u_n\|_{B_{p, r}^{s-1}}^2),\\
        \|(u_{0, n})^2\partial_x\rho_{0, n}\|_{B_{p, r}^{s-2}} &\lesssim\|\partial_x\rho_{0, n}\|_{B_{p, r}^{s-2}}\|(u_{0, n})^2\|_{B_{p, r}^{s-1}}\\
        &\lesssim \|\rho_{0, n}\|_{B_{p, r}^{s-1}}\|u_{0, n}\|_{L^\infty}\|u_{0, n}\|_{B_{p, r}^{s-1}}\lesssim 2^{-n(s+1)},\\
      \|\rho_{0, n}u_{0, n}\partial_xu_{0, n}\|_{B_{p, r}^{s-2}}&\lesssim\|\rho_{0, n}\|_{B_{p, r}^{s-2}}(\|u_{0, n}\|_{L^\infty}\|\partial_xu_{0, n}\|_{B_{p, r}^{s-1}}+\|u_{0, n}\|_{B_{p, r}^{s-1}}\|\partial_xu_{0, n}\|_{L^\infty})\\
      &\lesssim 2^{-n}(2^{-ns}+2^{-n}2^n2^{-ns}) \lesssim2^{-n(s+1)},\\
      \|[(u_n)^2-(u_{0, n})^2]\partial_xu_{0, n}\|_{B_{p, r}^{s-1}}
      &\lesssim \|u_n-u_{0, n}\|_{B_{p, r}^{s-1}}\|u_n+u_{0, n}\|_{B_{p, r}^{s-1}}
        \|\partial_xu_{0, n}\|_{B_{p, r}^{s-1}} \\
        &\lesssim \|u_n-u_{0, n}\|_{B_{p, r}^{s-1}}\|u_n, u_{0, n}\|^2_{B_{p, r}^s},\\
      \|[\mathcal{P}(u_n)-\mathcal{P}(u_{0, n})]\|_{B_{p, r}^{s-1}}&\lesssim\|u_n-u_{0, n}\|_{B_{p, r}^{s-1}}\|u_n, u_{0, n}\|^2_{B_{p, r}^s},\\
       \|[\mathcal{R}(u_n, \rho_n)-\mathcal{R}(u_{0, n}, \rho_{0, n})]\|_{B_{p, r}^{s-1}}&\lesssim\|u_n-u_{0, n}\|_{B_{p, r}^{s-1}}\|\rho_n\|^2_{B_{p, r}^{s-1}}\\
       &\;\;+\|\rho_n-\rho_{0, n}\|_{B_{p, r}^{s-2}}(\|u_{0, n}\|^2_{B_{p, r}^s}
       +\|\rho_n, \rho_{0, n}\|^2_{B_{p, r}^{s-1}}).
       \end{align*}
       Again using Lemma \ref{lem2.1} and the Banach algebra property of $B_{p, r}^{s-1}$, one has
       \begin{align*}
     \|\mathcal{P}_1(u_{0, n})\|_{B_{p, r}^{s-1}}&\lesssim\|(u_{0, n})^3\|_{B_{p, r}^{s-2}}\lesssim\|u_{0, n}\|_{B_{p, r}^{s-2}}\|(u_{0, n})^2\|_{B_{p, r}^{s-1}}\\
     &\lesssim \|u_{0, n}\|_{B_{p, r}^{s-2}}\|u_{0, n}\|_{L^\infty}\|u_{0, n}\|_{B_{p, r}^{s-1}}\lesssim 2^{-n(s+3)},\\
     \|\mathcal{P}_2(u_{0, n})\|_{B_{p, r}^{s-1}}&\lesssim\|u_{0, n}(\partial_xu_{0, n})^2\|_{B_{p, r}^{s-2}}\lesssim\|u_{0, n}\|_{B_{p, r}^{s-2}}\|(\partial_xu_{0, n})^2\|_{B_{p, r}^{s-1}}\\
     &\lesssim \|u_{0, n}\|_{B_{p, r}^{s-2}}\|\partial_xu_{0, n}\|_{L^\infty}\|\partial_xu_{0, n}\|_{B_{p, r}^{s-1}}\lesssim 2^{-n(s+1)},\\
     \|\mathcal{P}_3(u_{0, n})\|_{B_{p, r}^{s-1}}&\lesssim\|(\partial_xu_{0, n})^3\|_{B_{p, r}^{s-2}}\lesssim\|\partial_xu_{0, n}\|_{B_{p, r}^{s-2}}\|(\partial_xu_{0, n})^2\|_{B_{p, r}^{s-1}}\\
     &\lesssim \|u_{0, n}\|_{B_{p, r}^{s-1}}\|\partial_xu_{0, n}\|_{L^\infty}\|\partial_xu_{0, n}\|_{B_{p, r}^{s-1}}\lesssim 2^{-ns},\\
     \|\mathcal{R}_1(u_{0, n}, \rho_{0, n})\|_{B_{p, r}^{s-1}}&\lesssim\|u_{0, n}(\rho_{0, n})^2\|_{B_{p, r}^{s-2}}\lesssim \|u_{0, n}\|_{B_{p, r}^{s-2}}\|\rho_{0, n}\|_{L^\infty}\|\rho_{0, n}\|_{B_{p, r}^{s-1}}\lesssim 2^{-n(s+1)},\\
     \|\mathcal{R}_2(u_{0, n}, \rho_{0, n})\|_{B_{p, r}^{s-1}}&\lesssim\|\partial_xu_{0, n}(\rho_{0, n})^2\|_{B_{p, r}^{s-2}}\lesssim \|u_{0, n}\|_{B_{p, r}^{s-1}}\|\rho_{0, n}\|_{L^\infty}\|\rho_{0, n}\|_{B_{p, r}^{s-1}}\lesssim 2^{-ns},\\
     \|(u_{0, n})^2\partial_xu_{0, n}\|_{B_{p, r}^{s-1}}
      &\lesssim \|\partial_xu_{0, n}\|_{B_{p, r}^{s-1}} \|u_{0, n}\|_{L^\infty} \|u_{0, n}\|_{B_{p, r}^{s-1}}\lesssim 2^{-n(s+1)}.
        \end{align*}
For the term
\begin{eqnarray*}
\rho_nu_n\partial_xu_n-\rho_{0, n}u_{0, n}\partial_xu_{0, n}=(\rho_n-\rho_{0, n})u_n\partial_xu_n
+ \rho_{0, n}(u_n-u_{0, n})\partial_xu_n+\rho_{0, n}u_{0, n}\partial_x(u_n-u_{0, n}).
\end{eqnarray*}
Following the same procedure of estimates as above, we find that
\begin{align*}
\|(\rho_n-\rho_{0, n})u_n\partial_xu_n\|_{B_{p, r}^{s-2}}&\lesssim\|\rho_n-\rho_{0, n}\|_{B_{p, r}^{s-2}}\|u_n\|_{B_{p, r}^{s-1}}\|u_n\|_{B_{p, r}^s}\lesssim\|\rho_n-\rho_{0, n}\|_{B_{p, r}^{s-2}}\|u_n\|^2_{B_{p, r}^s},\\
\|\rho_{0, n}(u_n-u_{0, n})\partial_xu_n\|_{B_{p, r}^{s-2}}&\lesssim\|u_n-u_{0, n}\|_{B_{p, r}^{s-1}}(\|\rho_{0, n}\|^2_{B_{p, r}^{s-1}}+\|u_n\|^2_{B_{p, r}^s}),\\
\|\rho_{0, n}u_{0, n}\partial_x(u_n-u_{0, n})\|_{B_{p, r}^{s-2}}&\lesssim\|\partial_x(u_n-u_{0, n})\|_{B_{p, r}^{s-2}}\|\rho_{0, n}u_{0, n}\|_{B_{p, r}^{s-1}}\lesssim\|u_n-u_{0, n}\|_{B_{p, r}^{s-1}}\|\rho_{0, n}, u_{0, n}\|^2_{B_{p, r}^{s-1}}.
\end{align*}
Denote
\begin{eqnarray*}
X_s=\|\varrho\|_{B_{p, r}^{s-1}}+\|\epsilon\|_{B_{p, r}^s},
\end{eqnarray*}
taking the above estimates into (\ref{eq3.7})-(\ref{eq3.8}), we get
\begin{eqnarray*}
X_{s-1}\lesssim \int_0^tX_{s-1}(\|u_n, u_{0, n}\|^2_{B_{p, r}^s}+\|\rho_n, \rho_{0, n}\|^2_{B_{p, r}^{s-1}})+t2^{-ns},
\end{eqnarray*}
since $\{\rho_n, u_n\}$ is bounded in $B_{p, r}^{s-1}\times B_{p, r}^s$,
which together with the Gronwall Lemma imply
$$ X_{s-1}\lesssim 2^{-ns}.$$
Combining with (\ref{eq3.3}) for $k=1$ and the interpolation inequality, we obtain that
 $$X_s\lesssim X_{s-1}^{\frac{1}{2}}X_s^{\frac{1}{2}}\lesssim 2^{-\frac{ns}{2}}2^{\frac{n}{2}}\lesssim 2^{-\frac{n}{2}(s-1)}.$$
Thus we have complete the proof of Proposition \ref{pro1}.

In order to obtain the non-uniformly continuous dependence property for the system (\ref{eq3}), we will show that for the constructed  initial data $(\tilde{\rho}_{0, n}, v_{0, n})$ with small perturbation, it can not approximate to the solution
$(\tilde{\rho}_n, v_n)$.
\begin{proposition}\label{pro2}
 Under the assumptions of Theorem \ref{the1.1}, we have
      \begin{eqnarray}\label{eq3.9}
      \|\tilde{\rho}_n-\tilde{\rho}_{0, n}-t\mathbf{w}_0^n\|_{B_{p, r}^{s-1}}+\|v_n-v_{0, n}-t\mathbf{v}_0^n\|_{B_{p, r}^s}\lesssim t^2+2^{-\frac n2},
        \end{eqnarray}
here, $\mathbf{w}_0^n=(v_{0, n})^2\partial_x\tilde{\rho}_{0, n}$, $\mathbf{v}_0^n=(v_{0, n})^2\partial_xv_{0, n}.$
\end{proposition}
\noindent{\bf Proof} \; Firstly, with the help of (\ref{eq3.1})-(\ref{eq3.2}), (\ref{eq3.4}) and making full use of the product estimates in Lemma \ref{lem2.1}, for $\sigma\geq -\frac 12$, we have
\begin{align}
  \|\mathbf{v}_0^n\|_{B_{p, r}^{s+\sigma}}&\lesssim \|v_{0, n}\|^2_{L^\infty}\|\partial_xv_{0, n}\|_{B_{p, r}^{s+\sigma}}+ \|v_{0, n}\|_{L^\infty}\|v_{0, n}\|_{B_{p, r}^{s+\sigma}}\|\partial_xv_{0, n}\|_{L^\infty}\nonumber\\
  &\lesssim (2^{-\frac n2})^22^{n(\sigma+1)}+2^{-\frac n2}2^{n\sigma}2^{-\frac n2}\lesssim 2^{n\sigma},\label{eq3.10}\\
   \|\mathbf{w}_0^n\|_{B_{p, r}^{s+\sigma}}&\lesssim \|v_{0, n}\|^2_{L^\infty}\|\partial_x\tilde{\rho}_{0, n}\|_{B_{p, r}^{s+\sigma}}+ \|v_{0, n}\|_{L^\infty}\|v_{0, n}\|_{B_{p, r}^{s+\sigma}}\|\partial_x\tilde{\rho}_{0, n}\|_{L^\infty}\nonumber\\
  &\lesssim (2^{-\frac n2})^22^{n(\sigma+2)}+2^{-\frac n2}2^{n\sigma}2^{-\frac n2}\lesssim 2^{n(\sigma+1)},\label{eq3.11}\\
\|\mathbf{v}_0^n\|_{B_{p, r}^{s-1}}&\lesssim \|v_{0, n}\|^2_{B_{p, r}^{s-1}}\|\partial_xv_{0, n}\|_{B_{p, r}^{s-1}}\lesssim 2^{-n},\label{eq3.12}\\
\|\mathbf{w}_0^n\|_{B_{p, r}^{s-1}}&\lesssim \|v_{0, n}\|^2_{B_{p, r}^{s-1}}\|\partial_x\tilde{\rho}_{0, n}\|_{B_{p, r}^{s-1}}\leq C,\label{eq3.13}\\
\|\mathbf{w}_0^n\|_{B_{p, r}^{s-2}}&\lesssim \|\partial_x\tilde{\rho}_{0, n}\|_{B_{p, r}^{s-2}}\|v_{0, n}\|^2_{B_{p, r}^{s-1}}\lesssim 2^{-n}.\label{eq3.14}
\end{align}

Let
\begin{eqnarray*}
        \left\{\begin{array}{ll}
        z_n=\tilde{\rho}_n-\tilde{\rho}_{0, n}-t\mathbf{w}_0^n,\\
         \omega_n=v_n-v_{0, n}-t\mathbf{v}_0^n,\end{array}\right.
        \end{eqnarray*}
then we can derive from (\ref{eq3}) that $(z_n, \omega_n)$ satisfies
\begin{eqnarray}\label{eq3.15}
        \left\{\begin{array}{ll}
        \partial_tz_n-(v_n)^2\partial_xz_n=t(v_n)^2\partial_x\mathbf{w}_0^n+t(v_n+v_{0, n})\mathbf{v}_0^n\partial_x\tilde{\rho}_{0, n}\\
        \qquad\qquad\qquad\qquad\;\; +\omega_n(v_n+v_{0, n})\partial_x\tilde{\rho}_{0, n}+\tilde{\rho}_nv_n\partial_xv_n,\\
         \partial_t\omega_n-(v_n)^2\partial_x\omega_n=t(v_n+v_{0, n})\mathbf{v}_0^n\partial_xv_{0, n}+\omega_n(v_n+v_{0, n})\partial_xv_{0, n}\\
         \qquad\qquad\qquad\qquad\;\;+t(v_n)^2\partial_x\mathbf{v}_0^n +\mathcal{R}(v_n, \tilde{\rho}_n)+\mathcal{P}(v_n).
         \end{array}\right.
        \end{eqnarray}
Applying  Lemma \ref{lem2.1}, using (\ref{eq3.1})-(\ref{eq3.2}), (\ref{eq3.4}), (\ref{eq3.11})-(\ref{eq3.13}), we arrive at
\begin{align}
  \|(v_n)^2\partial_x\mathbf{w}_0^n\|_{B_{p, r}^{s-2}}&\lesssim \|\partial_x\mathbf{w}_0^n\|_{B_{p, r}^{s-2}}\|(v_n)^2\|_{B_{p, r}^{s-1}}\lesssim \|\mathbf{w}_0^n\|_{B_{p, r}^{s-1}}\|v_n\|^2_{B_{p, r}^{s-1}}\lesssim 2^{-n}, \label{eq3.16}\\
 \|(v_n)^2\partial_x\mathbf{w}_0^n\|_{B_{p, r}^{s-1}}&\lesssim \|(v_n)^2\|_{B_{p, r}^{s-1}}\|\partial_x\mathbf{w}_0^n\|_{B_{p, r}^{s-1}}\lesssim (2^{-\frac{n}{2}})^22^n\leq C, \label{eq3.17}\\
  \|(v_n+v_{0, n})\mathbf{v}_0^n\partial_x\tilde{\rho}_{0, n}\|_{B_{p, r}^{s-1}}&\lesssim \|v_n+v_{0, n}\|_{B_{p, r}^{s-1}}\|\mathbf{v}_0^n\|_{B_{p, r}^{s-1}}\|\partial_x\tilde{\rho}_{0, n}\|_{B_{p, r}^{s-1}}\lesssim 2^{-\frac{n}{2}}2^{-n}2^n\leq 2^{-\frac{n}{2}}, \label{eq3.18}\\
  \|\omega_n(v_n+v_{0, n})\partial_x\tilde{\rho}_{0, n}\|_{B_{p, r}^{s-2}}&\lesssim \|\partial_x\tilde{\rho}_{0, n}\|_{B_{p, r}^{s-2}}\|\omega_n\|_{B_{p, r}^{s-1}}\|v_n+v_{0, n}\|_{B_{p, r}^{s-1}}\lesssim \|\omega_n\|_{B_{p, r}^{s-1}}, \label{eq3.19}\\
  \|\omega_n(v_n+v_{0, n})\partial_x\tilde{\rho}_{0, n}\|_{B_{p, r}^{s-1}}&\lesssim \|\omega_n\|_{B_{p, r}^{s-1}}\|v_n+v_{0, n}\|_{B_{p, r}^{s-1}}\|\partial_x\tilde{\rho}_{0, n}\|_{B_{p, r}^{s-1}}\lesssim 2^{\frac{n}{2}}\|\omega_n\|_{B_{p, r}^{s-1}}, \label{eq3.20}\\
\|\tilde{\rho}_nv_n\partial_xv_n\|_{B_{p, r}^{s-2}}&\lesssim \|\tilde{\rho}_n\|_{B_{p, r}^{s-2}}\|v_n\|_{B_{p, r}^{s-1}}\|\partial_xv_n\|_{B_{p, r}^{s-1}}\lesssim 2^{-n}, \label{eq3.21}\\
\|\tilde{\rho}_nv_n\partial_xv_n\|_{B_{p, r}^{s-1}}&\lesssim \|\tilde{\rho}_n\|_{B_{p, r}^{s-1}}\|v_n\|_{B_{p, r}^{s-1}}\|\partial_xv_n\|_{B_{p, r}^{s-1}}\lesssim 2^{-\frac{n}{2}}, \label{eq3.22}
\end{align}
According to Lemma \ref{lem2.2} to (\ref{eq3.15}), using the fact that $\{v_n\}$ is bounded in $L^\infty_T(B_{p, r}^s)$, firstly with (\ref{eq3.16}), (\ref{eq3.18}), (\ref{eq3.19}), (\ref{eq3.21}), we infer that
\begin{equation}\label{eq3.23}
 \|z_n\|_{B_{p, r}^{s-2}}\leq C\int_0^t(\|z_n\|_{B_{p, r}^{s-2}}+\|\omega_n\|_{B_{p, r}^{s-1}})d\tau+Ct^22^{-\frac n2}+C2^{-n},
\end{equation}
and again combining with (\ref{eq3.17}), (\ref{eq3.18}), (\ref{eq3.20}), (\ref{eq3.22}), we obtain that
\begin{equation}\label{eq3.24}
 \|z_n\|_{B_{p, r}^{s-1}}\leq C\int_0^t\|z_n\|_{B_{p, r}^{s-1}}d\tau+C\int_0^t2^{\frac n2}\|\omega_n\|_{B_{p, r}^{s-1}}d\tau+Ct^2+C2^{-\frac{n}{2}}.
\end{equation}

In the following, we shall estimate $\omega_n$ in $B_{p, r}^{s-1}$ and $B_{p, r}^s$, respectively.
With the aid of Lemma \ref{lem2.1} and (\ref{eq3.2}), (\ref{eq3.4}), (\ref{eq3.10}), (\ref{eq3.12}), one has
\begin{align}
  \|(v_n+v_{0, n})\mathbf{v}_0^n\partial_xv_{0, n}\|_{B_{p, r}^s}&\lesssim \|v_n+v_{0, n}\|_{L^\infty}\|\mathbf{v}_0^n\partial_xv_{0, n}\|_{B_{p, r}^s}+\|v_n+v_{0, n}\|_{B_{p, r}^s}\|\mathbf{v}_0^n\partial_xv_{0, n}\|_{L^\infty}\nonumber\\
  &\lesssim \|v_n+v_{0, n}\|_{B_{p, r}^{s-1}}(\|\mathbf{v}_0^n\|_{L^\infty}\|v_{0, n}\|_{B_{p, r}^{s+1}}+\|\mathbf{v}_0^n\|_{B_{p, r}^s}\|\partial_xv_{0, n}\|_{L^\infty})\nonumber\\
&\;\;\; +\|v_n+v_{0, n}\|_{B_{p, r}^s}\|\mathbf{v}_0^n\|_{L^\infty}\|\partial_xv_{0, n}\|_{L^\infty}\nonumber\\
&\lesssim 2^{-\frac n2}(2^{-n}2^n+2^{-\frac n2})+2^{-n}2^{-\frac n2}\lesssim 2^{-\frac n2}, \label{eq3.25}\\
\|\omega_n(v_n+v_{0, n})\partial_xv_{0, n}\|_{B_{p, r}^s}&\lesssim \|\omega_n\|_{L^\infty}\|(v_n+v_{0, n})\partial_xv_{0, n}\|_{B_{p, r}^s}+\|\omega_n\|_{B_{p, r}^s}\|(v_n+v_{0, n})\partial_xv_{0, n}\|_{L^\infty}\nonumber\\
&\lesssim\|\omega_n\|_{B_{p, r}^{s-1}}(\|v_n+v_{0, n}\|_{L^\infty}\|v_{0, n}\|_{B_{p, r}^{s+1}}+\|v_n+v_{0, n}\|_{B_{p, r}^s}\|\partial_xv_{0, n}\|_{L^\infty})\nonumber\\
&\;\;\;+\|\omega_n\|_{B_{p, r}^s}\|v_n+v_{0, n}\|_{L^\infty}\|\partial_xv_{0, n}\|_{L^\infty}\nonumber\\
&\lesssim \|\omega_n\|_{B_{p, r}^{s-1}}(2^{-\frac n2}2^n+2^{-\frac n2})+\|\omega_n\|_{B_{p, r}^s}2^{-\frac n2}2^{-\frac n2}\nonumber\\
&\lesssim 2^{\frac n2}\|\omega_n\|_{B_{p, r}^{s-1}}+\|\omega_n\|_{B_{p, r}^s},\label{eq3.26}\\
\|\omega_n(v_n+v_{0, n})\partial_xv_{0, n}\|_{B_{p, r}^{s-1}}&\lesssim\|\omega_n\|_{B_{p, r}^{s-1}}\|v_n+v_{0, n}\|_{B_{p, r}^{s-1}}\|\partial_xv_{0, n}\|_{B_{p, r}^{s-1}}\lesssim\|\omega_n\|_{B_{p, r}^{s-1}}\label{eq3.27}\\
\|(v_n)^2\partial_x\mathbf{v}_0^n\|_{B_{p, r}^s}&\lesssim \|v_n\|^2_{B_{p, r}^{s-1}}\|\partial_x\mathbf{v}_0^n\|_{B_{p, r}^s}+\|v_n\|^2_{B_{p, r}^s}\|\partial_x\mathbf{v}_0^n\|_{L^\infty}\leq C, \label{eq3.28}\\
\|(v_n)^2\partial_x\mathbf{v}_0^n\|_{B_{p, r}^{s-1}}&\lesssim \|v_n\|^2_{B_{p, r}^{s-1}}\|\partial_x\mathbf{v}_0^n\|_{B_{p, r}^{s-1}}\lesssim 2^{-n}.\label{eq3.29}
\end{align}
For the term $\mathcal{R}(v_n, \tilde{\rho}_n)=\mathcal{R}_1(v_n, \tilde{\rho}_n)+\mathcal{R}_2(v_n, \tilde{\rho}_n),$ we have
\begin{align}
 \|\mathcal{R}_1(v_n, \tilde{\rho}_n)\|_{B_{p, r}^{s-1}}&\lesssim \|v_n(\tilde{\rho}_n)^2\|_{B_{p, r}^{s-2}}
\lesssim \|\tilde{\rho}_n\|_{B_{p, r}^{s-2}}\|\tilde{\rho}_n\|_{B_{p, r}^{s-1}}\|v_n\|_{B_{p, r}^{s-1}}\lesssim 2^{-n},\label{eq3.30}\\
 \|\mathcal{R}_1(v_n, \tilde{\rho}_n)\|_{B_{p, r}^s}&\lesssim \|v_n(\tilde{\rho}_n)^2\|_{B_{p, r}^{s-1}}
\lesssim \|\tilde{\rho}_n\|^2_{B_{p, r}^{s-1}}\|v_n\|_{B_{p, r}^{s-1}}\lesssim 2^{-\frac n2}.\label{eq3.31}
\end{align}
\begin{align*}
\mathcal{R}_2(v_n, \tilde{\rho}_n)=&-\underbrace{\frac 12(1-\partial_x^2)^{-1}(z_n(\tilde{\rho}_n+\tilde{\rho}_{0, n})\partial_xv_n)}_{\mathcal{R}_{2,1}}
-\underbrace{\frac 12(1-\partial_x^2)^{-1}(t\mathbf{w}_0^n(\tilde{\rho}_n+\tilde{\rho}_{0, n})\partial_xv_n)}_{\mathcal{R}_{2,2}}\\
&-\underbrace{\frac 12(1-\partial_x^2)^{-1}((\tilde{\rho}_{0, n})^2\partial_xv_n)}_{\mathcal{R}_{2,3}}.
\end{align*}
Using Lemma \ref{lem2.1} together with (\ref{eq3.1})-(\ref{eq3.2}), (\ref{eq3.4}), (\ref{eq3.10}), (\ref{eq3.14}), we have
\begin{align}
  \|\mathcal{R}_{2,1}\|_{B_{p, r}^s}\lesssim \|z_n(\tilde{\rho}_n+\tilde{\rho}_{0, n})\partial_xv_n\|_{B_{p, r}^{s-2}}\lesssim \|z_n\|_{B_{p, r}^{s-2}}\|\tilde{\rho}_n+\tilde{\rho}_{0, n}\|_{B_{p, r}^{s-1}}\|\partial_xv_n\|_{B_{p, r}^{s-1}}\lesssim \|z_n\|_{B_{p, r}^{s-2}} \label{eq3.32}\\
   \|\mathcal{R}_{2,2}\|_{B_{p, r}^s}\lesssim t\|\mathbf{w}_0^n(\tilde{\rho}_n+\tilde{\rho}_{0, n})\partial_xv_n\|_{B_{p, r}^{s-2}}\lesssim t\|\mathbf{w}_0^n\|_{B_{p, r}^{s-2}}\|\tilde{\rho}_n+\tilde{\rho}_{0, n}\|_{B_{p, r}^{s-1}}\|\partial_xv_n\|_{B_{p, r}^{s-1}}\lesssim t2^{-n}. \label{eq3.33}
\end{align}
Using the following estimate
\begin{align*}
 \|(\tilde{\rho}_{0, n})^2\partial_xv_{0, n}\|_{B_{p, r}^{s-2}}&\lesssim  \|(\tilde{\rho}_{0, n})^2\partial_xf_n\|_{B_{p, r}^{s-2}}+ \|(\tilde{\rho}_{0, n})^2\partial_xg_n\|_{B_{p, r}^{s-2}}\\
&\lesssim \|\tilde{\rho}_{0, n}\|_{B_{p, r}^{s-2}}(\|\tilde{\rho}_{0, n}\|_{L^\infty}\|\partial_xf_n\|_{B_{p, r}^{s-1}}+\|\tilde{\rho}_{0, n}\|_{B_{p, r}^{s-1}}\|\partial_xf_n\|_{L^\infty})\\
&\;\;+\|\tilde{\rho}_{0, n}\|_{B_{p, r}^{s-2}}(\|\tilde{\rho}_{0, n}\|_{L^\infty}\|\partial_xg_n\|_{B_{p, r}^{s-1}}+\|\tilde{\rho}_{0, n}\|_{B_{p, r}^{s-1}}\|\partial_xg_n\|_{L^\infty})\\
&\lesssim 2^{-\frac n2}(2^{-\frac n2}+2^{-ns}2^n)+2^{-\frac n2}(2^{-\frac n2}2^{-\frac n2}+2^{-\frac n2})\\
&\lesssim 2^{-n},
\end{align*}
hence we have
\begin{align}\label{eq3.34}
 \|\mathcal{R}_{2, 3}\|_{B_{p, r}^s}&\lesssim \|(\tilde{\rho}_{0, n})^2\partial_xv_n\|_{B_{p, r}^{s-2}}\nonumber\\
 &\lesssim \|(\tilde{\rho}_{0, n})^2\partial_xw_n\|_{B_{p, r}^{s-2}}+\|(\tilde{\rho}_{0, n})^2\partial_xv_{0,n}\|_{B_{p, r}^{s-2}}+\|t(\tilde{\rho}_{0, n})^2\partial_x\mathbf{v}_0^n\|_{B_{p, r}^{s-2}} \nonumber\\
 &\lesssim \|\omega_n\|_{B_{p, r}^{s-1}}+2^{-n}+t2^{-\frac n2}.
\end{align}
For the term $\mathcal{P}(v_n)=\mathcal{P}_1(v_n)+\mathcal{P}_2(v_n)+\mathcal{P}_3(v_n),$  we have from Lemma \ref{lem2.1} that
  \begin{equation}\label{eq3.35}
    \|\mathcal{P}_1(v_n)\|_{B_{p, r}^s}\lesssim \|(v_n)^3\|_{B_{p, r}^{s-1}}\lesssim 2^{-\frac {3}{2}n},
  \end{equation}
while it needs to be careful to deal with $\mathcal{P}_2(v_n)$ and $\mathcal{P}_3(v_n).$ By making full use of the structure of $v_n$, we find that
\begin{align*}
\mathcal{P}_2(v_n)=&\underbrace{\frac 32\partial_x(1-\partial_x^2)^{-1}(v_n\partial_x(v_n+v_{0, n})\partial_xw_n)}_{\mathcal{P}_{2, 1}}\\
&\;+\underbrace{\frac 32\partial_x(1-\partial_x^2)^{-1}(tv_n\partial_x(v_n+v_{0, n})\partial_x\mathbf{v}_0^n)}_{\mathcal{P}_{2, 2}}\\
&\;+\underbrace{\frac 32\partial_x(1-\partial_x^2)^{-1}(v_n(\partial_xv_{0, n})^2)}_{\mathcal{P}_{2, 3}},
\end{align*}
and
\begin{align}
  \|\mathcal{P}_{2, 1}\|_{B_{p, r}^s}&\lesssim \|v_n\partial_x(v_n+v_{0, n})\partial_xw_n\|_{B_{p, r}^{s-1}}\lesssim \|w_n\|_{B_{p, r}^s}, \label{eq3.36}\\
  \|\mathcal{P}_{2, 1}\|_{B_{p, r}^{s-1}}&\lesssim \|v_n\partial_x(v_n+v_{0, n})\partial_xw_n\|_{B_{p, r}^{s-2}}\lesssim \|w_n\|_{B_{p, r}^{s-1}}, \label{eq3.37}\\
 \|\mathcal{P}_{2, 2}\|_{B_{p, r}^s}&\lesssim t\|v_n\partial_x(v_n+v_{0, n})\partial_x\mathbf{v}_0^n\|_{B_{p, r}^{s-1}}\lesssim t2^{-\frac n2},\label{eq3.38}
\end{align}
\begin{align}
  \|\mathcal{P}_{2, 3}\|_{B_{p, r}^s}&\lesssim \|v_n(\partial_xv_{0, n})^2\|_{_{B_{p, r}^{s-1}}}\nonumber\\
  &\lesssim \|v_n\|_{L^\infty}\|\partial_xv_{0, n}\|_{L^\infty}\|\partial_xv_{0, n}\|_{_{B_{p, r}^{s-1}}}+
  \|v_n\|_{B_{p, r}^{s-1}}\|\partial_xv_{0, n}\|^2_{L^\infty}\nonumber\\
  &\lesssim 2^{-\frac n2}(2^{-ns}2^n+2^{-\frac n2})+2^{-\frac n2}(2^{-ns}2^n+2^{-\frac n2})^2\nonumber\\
  &\lesssim 2^{-n}. \label{eq3.39}
\end{align}
$\mathcal{P}_3(v_n)$ can be performed in a similar way. Firstly, we obtain that
\begin{align*}
\mathcal{P}_3(v_n)=&\underbrace{\frac 12(1-\partial_x^2)^{-1}(\partial_xv_n\partial_x(v_n+v_{0, n})\partial_xw_n)}_{\mathcal{P}_{3, 1}}\\
&\;+\underbrace{\frac 12(1-\partial_x^2)^{-1}(t\partial_xv_n\partial_x(v_n+v_{0, n})\partial_x\mathbf{v}_0^n)}_{\mathcal{P}_{3, 2}}\\
&\;+\underbrace{\frac 12(1-\partial_x^2)^{-1}(\partial_xv_n(\partial_xv_{0, n})^2)}_{\mathcal{P}_{3, 3}},
\end{align*}
then we have from Lemma \ref{lem2.1} together with (\ref{eq3.4}) and (\ref{eq3.12}) that
\begin{align}
  \|\mathcal{P}_{3, 1}\|_{B_{p, r}^s}&\lesssim \|\partial_xv_n\partial_x(v_n+v_{0, n})\partial_xw_n\|_{B_{p, r}^{s-2}}\lesssim \|w_n\|_{B_{p, r}^{s-1}}, \label{eq3.40}\\
  \|\mathcal{P}_{3, 2}\|_{B_{p, r}^s}&\lesssim t\|\partial_xv_n\partial_x(v_n+v_{0, n})\partial_x\mathbf{v}_0^n\|_{B_{p, r}^{s-2}}\nonumber\\
  &\lesssim t\|\partial_x\mathbf{v}_0^n\|_{B_{p, r}^{s-2}}\|\partial_xv_n\|_{B_{p, r}^{s-1}}\|\partial_x(v_n+v_{0, n})\|_{B_{p, r}^{s-1}}\lesssim t2^{-n}, \label{eq3.41}\\
\|\mathcal{P}_{3, 3}\|_{B_{p, r}^s}&\lesssim \|\partial_xv_n(\partial_xv_{0, n})^2\|_{B_{p, r}^{s-2}}\lesssim \|\partial_xv_n\|_{B_{p, r}^{s-2}}\|\partial_xv_{0, n}\|_{L^\infty}\|\partial_xv_{0, n}\|_{B_{p, r}^{s-1}}\lesssim 2^{-n}. \label{eq3.42}
\end{align}
Applying Lemma \ref{lem2.2} firstly together with (\ref{eq3.25}), (\ref{eq3.27}), (\ref{eq3.29}),  (\ref{eq3.30}),
(\ref{eq3.32})-(\ref{eq3.35}), (\ref{eq3.37})-(\ref{eq3.42})  to (\ref{eq3.15}), using the fact that $\{v_n\}$ is bounded in $L^\infty_T(B_{p, r}^s)$, we infer that
\begin{equation}\label{eq3.43}
 \|\omega_n\|_{B_{p, r}^{s-1}}\leq C\int_0^t(\|z_n\|_{B_{p, r}^{s-2}}+\|\omega_n\|_{B_{p, r}^{s-1}})d\tau+Ct^22^{-\frac n2}+C2^{-n},
\end{equation}
and again combining with (\ref{eq3.25}), (\ref{eq3.26}), (\ref{eq3.28}),
(\ref{eq3.31})-(\ref{eq3.36}), (\ref{eq3.38})-(\ref{eq3.42}), we obtain
\begin{equation}\label{eq3.44}
 \|\omega_n\|_{B_{p, r}^s}\leq C\int_0^t\|\omega_n\|_{B_{p, r}^s}d\tau+ C\int_0^t2^{\frac n2}(\|z_n\|_{B_{p, r}^{s-2}}+\|\omega_n\|_{B_{p, r}^{s-1}})d\tau+Ct^22^{-\frac n2}+C2^{-\frac n2}.
\end{equation}
Using Gronwall Lemma to (\ref{eq3.23}) and (\ref{eq3.43}) imply
\begin{align*}
 \|z_n\|_{B_{p, r}^{s-2}}+\|\omega_n\|_{B_{p, r}^{s-1}}\leq Ct^22^{-\frac n2}+C2^{-n},
\end{align*}
which together with (\ref{eq3.24}) and (\ref{eq3.44}) yield that
 \begin{equation*}
  \|z_n\|_{B_{p, r}^{s-1}}+\|\omega_n\|_{B_{p, r}^s}\leq Ct^2+C2^{-\frac n2}.
\end{equation*}
 Thus, we have finished the proof of Proposition \ref{pro2}.

{\bf Proof of Theorem \ref{the1.1}} It is obvious that
 \begin{align*}
   \|\tilde{\rho}_{0, n}-\rho_{0, n}\|_{B_{p, r}^{s-1}}&=\|g_n\|_{B_{p, r}^{s-1}}\leq C2^{-\frac n2},\\
   \|v_{0, n}-u_{0, n}\|_{B_{p, r}^s}&=\|g_n\|_{B_{p, r}^s}\leq C2^{-\frac n2},
 \end{align*}
 which mean that
 \begin{equation*}
   \lim_{n\rightarrow \infty}(\|\tilde{\rho}_{0, n}-\rho_{0, n}\|_{B_{p, r}^{s-1}}+\|v_{0, n}-u_{0, n}\|_{B_{p, r}^s})=0.
 \end{equation*}
However, according to Proposition \ref{pro1} and Proposition \ref{pro2}, we get
\begin{align}\label{eq3.45}
  \|\tilde{\rho}_n-\rho_n\|_{B_{p, r}^{s-1}}&=\|z_n+t\mathbf{w}_0^n+g_n+\rho_{0, n}-\rho_n\|_{B_{p, r}^{s-1}}\nonumber\\
  &\gtrsim t\|\mathbf{w}_0^n\|_{B_{p, r}^{s-1}}-t^2-2^{-\frac n2}-2^{-\frac n2 (s-1)}.
\end{align}
Notice that
\begin{align*}
 \mathbf{w}_0^n&=(v_{0, n})^2\partial_x\tilde{\rho}_{0, n}\\
 &=(f_n)^2\partial_x(2^nf_n)+g_n^2\partial_x(2^nf_n)+(2f_ng_n)\partial_x(2^nf_n)+(v_{0, n})^2\partial_xg_n.
\end{align*}
With the aid of Lemma \ref{lem2.1} and the Banach algebra property of $B_{p, r}^{s-1}$, we find that
\begin{align*}
 \|(f_n)^2\partial_x(2^nf_n)\|_{B_{p, r}^{s-1}}&\lesssim \|f_n\|^2_{B_{p, r}^{s-1}}\|\partial_x(2^nf_n)\|_{B_{p, r}^{s-1}}\lesssim (2^{-n})^22^n\lesssim 2^{-n},\\
\|(2f_ng_n)\partial_x(2^nf_n)\|_{B_{p, r}^{s-1}}&\lesssim \|f_n\|_{B_{p, r}^{s-1}}\|g_n\|_{B_{p, r}^{s-1}}\|2^nf_n\|_{B_{p, r}^s}\lesssim 2^{-n}2^{-\frac n2}2^n\lesssim 2^{-\frac n2},\\
\|(v_{0, n})^2\partial_xg_n\|_{B_{p, r}^{s-1}}&\lesssim \|v_{0, n}\|^2_{B_{p, r}^{s-1}}\|\partial_xg_n\|_{B_{p, r}^{s-1}}\lesssim 2^{-\frac 32 n}.
\end{align*}
However, using the fact that $\Delta_j\big((g_n)^2\partial_x(2^nf_n)\big)=0, j\neq n$ and $\Delta_n\big((g_n)^2\partial_x(2^nf_n)\big)=(g_n)^2\partial_x(2^nf_n)$ for $n\geq 5$, direct calculation
 shows that
\begin{align*}
  &\|g_n^2\partial_x(2^nf_n)\|_{B_{p, r}^{s-1}}=2^{n(s-1)}\|g_n^2\partial_x(2^nf_n)\|_{L^p}\\
  =&\|2^{-n}\phi^2(x)\partial_x\phi(x)\sin(\frac{17}{12}2^nx)+\frac{17}{12}\phi^3(x)\cos(\frac{17}{12}2^nx)\|_{L^p}\\
  \gtrsim &\|\frac{17}{12}\phi^3(x)\cos(\frac{17}{12}2^nx)\|_{L^p}-2^{-n}\rightarrow \frac{17}{12} \big(\frac{\int_0^{2\pi}|\cos x|^pdx}{2\pi}\big)^{\frac 1p}\|\phi^3(x)\|_{L^p},
\end{align*}
by the Riemann Theorem.
 
 Taking the above estimates into (\ref{eq3.45}) yields
  \begin{align*}
   \liminf_{n\rightarrow \infty}||\tilde{\rho}_n-\rho_n||_{B^{s-1}_{p,r}}\gtrsim t\quad\text{for} \ t \ \text{small enough}.
  \end{align*}
 Similarly, we have 
   \begin{align*}
   \liminf_{n\rightarrow \infty}||v_n-u_n||_{B^s_{p,r}}\gtrsim t\quad\text{for} \ t \ \text{small enough}.
  \end{align*}
  
  This completes the proof of Theorem \ref{the1.1}.

\section*{Acknowledgments}
 This work is partially supported by the National Natural Science Foundation of China (Grant No.11801090).

\end{document}